\documentclass[12pt,a4paper,twoside,leqno]{article}
\usepackage{graphicx,color}
\usepackage{amsthm,amsmath,amssymb}
\usepackage{overpic}
\usepackage{subfigure}
\newtheorem{teor}{Theorem}[section]
\newtheorem{defin}[teor]{Definition}
\newtheorem{lemm}[teor]{Lemma}
\newtheorem{osse}[teor]{Remark}
\newtheorem{prop}[teor]{Proposition}
\newtheorem{defi}[teor]{Definition}
\newtheorem{coro}[teor]{Corollary}
\newtheorem{prob}[teor]{Problem}
\newtheorem{hypo}[teor]{Hypothesis}

\newcommand{\bele}{\begin{lemm}\begin{sl}}
\newcommand{\enle}{\end{sl}\end{lemm}}
\newcommand{\bedef}{\begin{defi}\begin{sl}}
\newcommand{\eddef}{\end{sl}\end{defi}}
\newcommand{\bete}{\begin{teor}\begin{sl}}
\newcommand{\ente}{\end{sl}\end{teor}}
\newcommand{\beos}{\begin{osse}\begin{sl}}
\newcommand{\eddos}{\end{sl}\end{osse}}
\newcommand{\bepr}{\begin{prop}\begin{sl}}
\newcommand{\empr}{\end{sl}\end{prop}}
\newcommand{\bepro}{\begin{prob}\begin{rm}}
\newcommand{\empro}{\end{rm}\end{prob}}
\newcommand{\bede}{\begin{defin}\begin{sl}}
\newcommand{\edde}{\end{sl}\end{defin}}
\newcommand{\beco}{\begin{coro}\begin{sl}}
\newcommand{\enco}{\end{sl}\end{coro}}

\newcommand{\behy}{\begin{hypo}\begin{sl}}
\newcommand{\enhy}{\end{sl}\end{hypo}}

\def\dd{\;\!\mathrm{d}} % differential for integration

\newcommand{\RR}{\mathbb{R}}

\newcommand{\CC}{\mathbb{C}}
\newcommand{\LL}{\mathbb{L}}

\newcommand{\WW}{\mathbb{W}}
\newcommand{\HH}{\mathbb{H}}
\newcommand{\XX}{\mathbb{X}}
\newcommand{\OO}{\mathbb{O}}

\newcommand{\beeq}[1]{\begin{equation}\label{#1}}
\newcommand{\eddeq}{\end{equation}}

\newcommand{\beeqa}[1]{\begin{eqnarray}\label{#1}}
\newcommand{\eddeqa}{\end{eqnarray}}

\newcommand{\beal}[1]{\begin{align}\label{#1}}
\newcommand{\eddal}{\end{align}}

\newcommand{\bespl}[1]{\begin{split}\label{#1}}
\newcommand{\edspl}{\end{split}}

\newcommand{\bega}[1]{\begin{gather}\label{#1}}
\newcommand{\edga}{\end{gather}}

\newcommand{\beeqax}{\begin{eqnarray*}}
\newcommand{\eddeqax}{\end{eqnarray*}}

\def\qed{\ifmmode % if math mode, assume display: omit penalty etc.
 \else \leavevmode\unskip\penalty9999 \hbox{}\nobreak\hfill
 \fi
 \quad\hbox{\hskip.5em\vrule width.4em height.6em depth.05em\hskip.1em}}
\def\endproofsym{\qed}

\def\endnobox{\def\endproofsym{}\end{proof}\def\endproofsym{\qed}}

\newcommand{\no}{\nonumber}

\newcommand{\beeqao}{\begin{eqnarray}\no}
\newcommand{\bealo}{\begin{align}\no}
\newcommand{\besplo}{\begin{split}\no}
\newcommand{\begao}{\begin{gather}\no}

\newcommand{\duav}[1]{\langle{#1}\rangle}

\newcommand{\itt}{\int_0^t}

\DeclareMathOperator{\dive}{div}
\DeclareMathOperator{\curl}{curl}

\newcommand{\teta}{\vartheta}

\numberwithin{equation}{section}
\begin{document}

\title{Analysis and simulations\\
 of multifrequency induction hardening\footnote{A part of the research this work is based on has been funded by the German Federal Ministry of Education and Research in the framework of the program "Mathematik f\"ur Innovationen in Industrie und Dienstleistungen".}}

\author{
Dietmar H\"omberg \footnote{Weierstrass Institute for Applied
Analysis and Stochastics, Mohrenstr.~39, D-10117 Berlin, Germany,
E-mail {\tt  hoemberg@wias-berlin.de}. A part of this work was
done during D.~H\"omberg's visit at the Department of Mathematics
of the University of Milan in November 2011, supported
by the FP7-IDEAS-ERC-StG Grant \#256872 (EntroPhase).} \and 
Thomas Petzold\footnote{Weierstrass Institute for Applied
Analysis and Stochastics, Mohrenstr.~39, D-10117 Berlin, Germany,
E-mail {\tt  petzold@wias-berlin.de}. }
\and Elisabetta Rocca
\footnote{Weierstrass Institute for Applied
Analysis and Stochastics, Mohrenstr.~39, D-10117 Berlin, Germany,
E-mail {\tt  rocca@wias-berlin.de} and Dipartimento di Matematica, Universit\`a di Milano, Via
Saldini 50, 20133 Milano, Italy, E-mail {\tt
elisabetta.rocca@unimi.it}. The work of E.R. was supported
by the FP7-IDEAS-ERC-StG Grant \#256872 (EntroPhase). A part of this work was done
during E.~Rocca's visit at WIAS Berlin in May 2012.}}

\date{}

\maketitle

\begin{abstract}
We study a model for induction hardening of steel. The related differential system consists of
a time domain vector potential formulation of the Maxwell's equations coupled with an internal energy balance
and an ODE for the volume fraction of {\sl austenite}, the high temperature phase in steel. We first
solve the initial boundary value problem associated by means of a Schauder fixed point argument coupled with
suitable a-priori estimates and regularity results. Moreover, we prove a stability estimate entailing, in
particular, uniqueness of solutions for our Cauchy problem. We conclude with some finite element simulations for the coupled system.
\end{abstract}

%%\tableofcontents

%\newpage

%%%%%%%%%%%%%%%%%%%%%%%%%%%%%%%%%%%%%%%%%%%%%%%%%%%%%%%%%%%%%%%%%%%

%%%%%%%%%%%%%%%%%%%%%%%%%%%%%%%%%%%%%%%%%%%%%%%%%%%%%%%%%%%%%%%%%%%

%%%%%%%%%%%%%%%%%%%%%%%%%%%%%%%%%%%%%%%%%%%%%%%%%%%%%%%%%%%%%%%%%%%
%%%%%%%%%%%%%%%%%%%%%%%%% Introduction %%%%%%%%%%%%%%%%%%%%%%%%%%%%

\section{Introduction}
\label{intro}

In induction hardening a coil that is connected to an alternating current source generates a periodically changing electro-magnetic field. The temporal changing magnetic flux induces a current in the workpiece that is enclosed by the induction coil. Due to the resistance of the workpiece, some part of the power is transformed into eddy current losses that result in Joule heating. 
The latter leads to a change of microstructure to the high temperature phase in steel called austenite. After switching off the current and possibly a short holding time the workpiece is cooled rapidly and the austenite layer produced upon heating is transformed to another phase called martensite, responsible for the desired hardening effect.

Induction heat treatments can easily be integrated into a process chain. Moreover, they are energy efficient since the heat is generated directly in the workpiece. That is why induction hardening is still the most important surface  treatment technology.

Due to the skin effect, the eddy currents tend to distribute in a small surface layer. 
The penetration depth of these eddy currents depends on the material and essentially on the frequency. Therefore, it is difficult to obtain a uniform contour hardened zone for complex workpiece geometries such as gears using a current with only one frequency.
If for example, a high frequency (HF) is applied, then the penetration depth is small and it is possible to harden only the tip of the gear. With a medium frequency (MF) it is possible to heat the root of the gear, but not the tip. With a single frequency, a hardening of the complete tooth can only be achieved by increasing the heating time.  But then, the complete tooth is heated beyond the austenitization temperature, which results in a complete martensitic structure of the tooth after quenching, which is not desirable, since this will foster fatigue effects. 

\begin{figure}[htb]
\centering
\includegraphics[width=.7\linewidth]{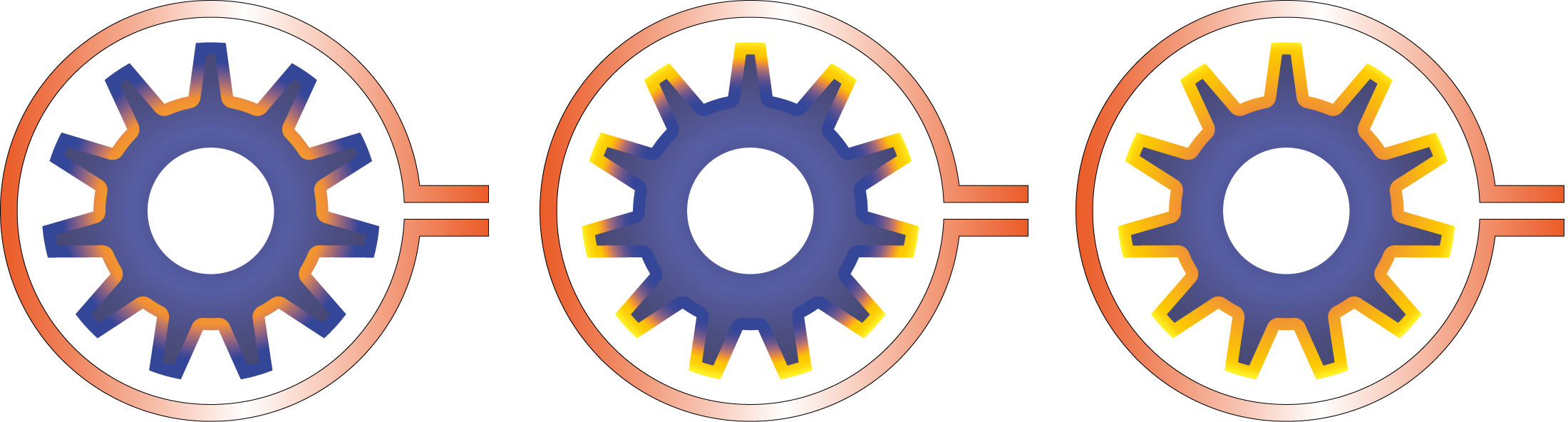}%
\caption{The effect of medium-, high- and multifrequency induction
  heating. MF (left): only the root of the gear is heated, HF
  (middle): only the tip of the gear is heated, MF+HF (right): tip and
  root of the gear are heated.} 
\label{fig:multifrequency}%
\end{figure}

Recently, a new approach has been developed which amounts to supplying medium and high frequency powers simultaneously on the induction coil. This concept is called \emph{multifrequency induction hardening}, see also Figure \ref{fig:multifrequency}.
%
%In order to achieve a hardening profile that follows the contour of the gear, very short heating times are necessary in order to avoid heat diffusion into the workpiece.
%In the past, multifrequency hardening was performed by subsequent heating of gearwheels in two separate inductors that are fed by power supplies with different frequencies. 

%Newer developments in dual-frequency induction hardening do not require a frequency changeover. MF and HF energy are supplied simultaneously to one inductor. 
The inductor current consists of a medium frequency fundamental oscillation superimposed by a high frequency oscillation. The amplitudes of both frequencies are independently controllable, which allows separate regulation of the respective shares of the output power of both frequencies according to the requirements of the workpiece. This provides the ability to control the depth of hardening at the root and the tip of the tooth  individually \cite{schwenk:2003}.

The main building blocks for a mathematical model of multifrequency hardening are an eddy current formulation of the Maxwell equations in the time domain, coupled to the balance of internal energy to describe the temperature evolution, and a model for the evolution of the high temperature phase austenite. Since the electric conductivity and the magnetic permeability may depend on temperature and/or the phase fractions and the Joule effect in the time domain is modelled  by the square of the time derivative of the magnetic vector potential we are faced with a strongly coupled nonlinear system of evolution equations. 

In two recent papers the simpler frequency domain situation of Joule heating has been studied. In \cite{GO14} the Boccardo-Galluet approach has been used to prove existence of a weak solution, while in \cite{DKSTY11}  new regularity results established in \cite{druet07} have been used to prove existence and stability in the frequency domain setting. 
In \cite{hoemberg04, HoembergSokolowski03} the eddy current model has been considered in the time domain. In the former the existence of a weak solution to the fully coupled model has been proven, in the latter also stability results are established, but only for a model with one-sided coupling in which the electric conductivity and permeability are assumed to be constant. 

The main novelty of the present paper is an existence and stability result for the strongly coupled time domain eddy current and Joule heating system. The paper is organised as follows. In the next section we derive the model equations and show that it complies with the second law of thermodynamics in form of the Clausius-Duhem inequality. In Section 3 we formulate the main mathematical results, the proofs of which are given in Section 4. The last section is devoted to presenting some results of numerical simulations for the coupled system based on a two time step approach with sequential decoupling of the evolution equations.

%%%%%%%%%%%%%%%%%%%%%%%%%%%%%%%%%%%%%%%%%%%%%%%%%%%%%%%%%%%%%%%%%%%
%%%%%%%%%%%%%%%%%%%%%%%%% model %%%%%%%%%%%%%%%%%%%%%%%%%%%%%%

\parskip0mm

\section{The model}
\label{model}

We restrict to the following idealized geometric setting (cf. the following Figure~\ref{domain}). Let $D\in \RR^3$ be a
domain containing the inductor coil $\Omega$ and the workpiece $\Sigma$. Assume
that $\overline{\Omega}\subset D$, $\overline{\Sigma}\subset D$, $\overline{\Omega}\cap\overline{\Sigma}=\emptyset$
and $\partial\Omega$, $\partial\Sigma$, $\partial D$ are of class $C^{1,1}$. Call $G=\Omega\cup\Sigma$ the set of conductors and
define the space-time domain as $Q=\Sigma\times (0,T)$.

%\begin{figure}[htbp] %  figure placement: here, top, bottom, or page
%  \centering
%  \includegraphics[scale=0.3]{gebiet.pdf}
%  \caption{The domain}
%  \label{domain}
%\end{figure}
\begin{figure}[htb]
  \centering
  \begin{overpic}[width=50mm]{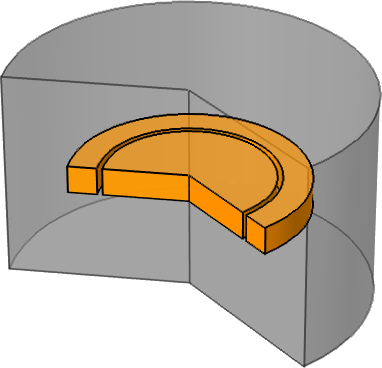}
  \put(50,50){$\Sigma$}
  \put(30,30){$D$}
  \put(70,40){$\Omega$}
  \end{overpic}
\caption{Domain $D$ consisting of the inductor $\Omega$, the workpiece
  $\Sigma$ and the surrounding air.} 
  \label{domain}
\end{figure}

Following the model derivation
of \cite{hoemberg04}, in the eddy current problems we neglect displacement currents and so we get
the following Maxwell's equations (cf. \cite{DL}) in $D\times (0,T)$:
\begin{align}
\label{M1}
&\curl H=J,\\
\label{M2}
&\curl E=-B_t,\\
\label{M3}
&\dive B=0\,,
\end{align}
where $E$ is the electric field, $B$ the magnetic induction, $H$ the magnetic field, $J$ the spatial current density.

We include phase transitions by describing the evolution of the volume fraction of austenite $z$
with the following problem derived in \cite{LeblondDevaux84} and \cite{FuhrmannHoemberg99}:
\begin{align}\label{eqz}
&z_t(t)= \frac{\left(z_{eq}(\teta)-z\right)^+}{\tau(\teta)}\quad \hbox{in } Q\,,\\
\no
&z(0)=0\quad \hbox{in }\Sigma\,.
\end{align}
We will  make precise the assumption on $\tau$ and $z_{eq}$ in the following section, however note that
$z_{eq}\in [0,1]$ is an equilibrium fraction
of {\sl austenite}, $\tau(\teta)$ a time constant, and by $(\cdot)^+$ we denote the positive part function.

Then, we assume the Ohm's law and a linear relation between the magnetic induction and the
magnetic field in $D\times (0,T)$:
\begin{align}\label{ohm}
& J =\sigma E\,,\\
\no
& B=\mu H\,,
\end{align}
where the electrical conductivity $\sigma$ and the magnetic permeability $\mu$ (sufficiently regular and bounded from below and above)
may depend both on the spatial variables and also
on the phase parameter $z$ and on the absolute temperature $\teta$:
\begin{equation}\label{def:sigma}
\sigma(x,z, \teta)=\begin{cases}
0,\quad& x\in D\setminus G\,,\\
\sigma_w(z, \teta),\quad& x\in \Sigma\,,\\
\sigma_i,\quad& x\in \Omega\,,
\end{cases}
\end{equation}
with its derivative with respect to the $z$-variable 
\begin{equation}\label{def:sigmaz}
\sigma_z(x,z, \teta)=\begin{cases}
0,\quad& x\in D\setminus \Sigma\,,\\
\sigma_w'(z, \teta),\quad& x\in \Sigma\,,
\end{cases}
\end{equation}
and
\begin{equation}\label{def:mu}
\mu(x,z, \teta)=\begin{cases}
\mu_0,\quad& x\in D\setminus G\,,\\
\mu_w(z, \teta),\quad& x\in \Sigma\,,\\
\mu_i,\quad& x\in \Omega\,,
\end{cases}
\end{equation}
with
\begin{equation}\label{def:muz}
\mu_z(x,z, \teta)=\begin{cases}
0,\quad& x\in D\setminus \Sigma\,,\\
\mu_w'(z, \teta),\quad& x\in \Sigma\,.
\end{cases}
\end{equation}
In view of \eqref{M3}, we introduce the magnetic vector potential $A$ such that
\begin{equation}\label{defA}
B=\curl A\quad \hbox{in }D
\end{equation}
and, since $A$ is not uniquely defined, we impose the Coulomb gauge
\[
\dive A=0 \quad \hbox{in }D\,.
\]
Using \eqref{M2} and \eqref{defA}, we define the scalar potential $\phi$ by
\begin{equation}\label{defPhi}
E+A_t=-\nabla\phi\quad\hbox{in }D\times (0,T)\,.
\end{equation}
Using \eqref{ohm}, we then obtain the following form for the total current density
\begin{equation}\label{totalcurrent}
J=-\sigma A_t-\sigma \nabla\phi\quad\hbox{in }D\times (0,T)\,,
\end{equation}
which, together with \eqref{defA} and \eqref{M1}, gives
\begin{equation}\label{eqA}
\sigma A_t+\curl\left(\frac{1}{\mu}\curl A\right)=J_{source} \quad \hbox{in }D\times (0,T)\,,
\end{equation}
where, for a given coil geometry (here a torus with rectangular cross-section),
the source current density
\begin{equation}\label{defu}
J_{source}=-\sigma \nabla \phi
\end{equation}
can be precomputed analytically
and it can be used as control for optimization, i.e., it can be taken of the form
$J_{source}=u(t)J_0$, where $J_0$ is the spatial current density prescribed in the induction  coil $\Omega$
\begin{equation}\label{def:J0}
J_0(x)=\begin{cases}
J_i(x), \quad &x\in \Omega\,\\
0,\quad&x\in D\setminus \Omega\,,
\end{cases}
\end{equation}
and $u=u(t)$ denotes a time-dependent control on $[0,T]$.
Assuming to have as constant density $\rho=1$ (for simplicity), the internal energy balance  in $Q$ results as
\begin{equation}\label{inten}
e_t+\dive q=JE= \sigma|A_t+\nabla\phi|^2=\sigma|A_t|^2\,,
\end{equation}
where $e$ denotes the internal energy of the system and $q$ the heat flux, which, accordingly to the standard
Fourier law is assumed as follows
\begin{equation}\label{fourier}
q=-\kappa(\teta)\nabla\teta, \quad \kappa>0\,.
\end{equation}
In case $\kappa(\teta)=\kappa>0$ constant, we get the standard Fourier law. 
From the Helmholtz relation $e=\psi+\teta s$, where $\psi=\psi(\teta, z)$ denotes the free energy of the system,
\eqref{inten} and \eqref{fourier},
we have that the Clausius-Duhem inequality
\begin{align}\label{clausius}
\teta \left(s_t+\dive \left(\frac{q}{\teta}\right)\right)&=e_t-\psi_\teta\teta_t-\psi_zz_t-\teta_t s-\frac{ q}{\teta}\nabla\teta+\dive q\\
\no
&=e_t+\dive q-(\psi_\teta+s)\teta_t-\psi_zz_t+\frac{\kappa|\nabla\teta|^2}{\teta}\\
\no
&=\sigma|A_t|^2-(\psi_\teta+s)\teta_t-\psi_zz_t+\frac{\kappa|\nabla\teta|^2}{\teta}\\
\no
&\geq 0
\end{align}
is satisfied e.g. if we assume the standard relations $\psi_\teta+s=0$ and $L\tau(\teta)z_t=-\psi_z$,
and hence (cf.\eqref{eqz}) $\psi_z= -L(z_{eq}(\teta)-z)^+$, for
some positive constant $L>0$.
Then, using the definition of the specific heat $c_v=\teta s_\teta$ (cf., e.g., \cite{LemaitreChaboche90}), we get
\[
e_t=c_v\teta_t+(\psi_z+\teta s_z)z_t=c_v\teta_t+\mathcal{F}(\teta,z)z_t\,,
\]
where we have denoted for simplicity  by
\begin{equation}\label{def:f}
\mathcal{F}(\teta, z)=\psi_z-\teta(\psi_z)_\teta=-L(z_{eq}(\teta)-z)^++L\teta z_{eq}'(\teta)\mathcal{H}(z_{eq}(\teta)-z)\,.
\end{equation}
Here $\mathcal{H}$ represents the Heaviside function.
Hence, the internal energy balance \eqref{inten} can be rewritten as
\begin{equation}\label{eqteta}
c_v\teta_t+\dive q=\sigma|A_t|^2-\mathcal{F}(\teta, z)z_t\quad \hbox{in }Q.
\end{equation}
%and $s=-\psi_\teta$ the entropy,
%we deduce $e_t=\teta s_t+\psi_z z_t=\teta s_\teta\teta_t+\teta s_zz_t+\psi_z z_t$.
%Let us assume that the latent heat $L=\psi_z$ is constant,
%then,  $s_z=-(\psi_z)_\teta=0$.
%Taking, moreover a constant specific heat (cf. \cite{LemaitreChaboche90}) $c_v=\teta s_\teta$,
%we get $e_t=c_v\teta_t+Lz_t$. Hence, the internal energy balance \eqref{inten} can be rewritten as
%\begin{equation}\label{eqteta}
%c_v\teta_t+Lz_t+\dive q=\sigma|A_t|^2\quad \hbox{in }Q.
%\end{equation}
%Let us note that this model turns out to satisfy the Clausius-Duhem inequality if we assume, e.g., $L<0$. Indeed in this case for all
%solutions of the field equations \eqref{eqz}, \eqref{eqA}, \eqref{inten}, \eqref{fourier}, the following
%inequality holds true
%\[
%\teta \left(s_t+\dive \left(\frac{q}{\teta}\right)\right)=e_t-Lz_t-\frac{ q}{\teta}\nabla\teta+\dive q=\sigma|A_t+\nabla\phi|^2-L f(\teta, z)
%+\kappa\frac{|\nabla\teta|^2}{\teta}\geq 0\,.
%\]
%More in general we could assume the following dynamic for $z$:  $\tau(\teta)z_t=-\psi_z$ in order to preserve the
%thermodynamical consistency of the model, but we do not want to enter this topic here.
Finally, we couple the system (\ref{eqz}), \eqref{eqA}, \eqref{defu}, (\ref{fourier}), (\ref{def:f}--\ref{eqteta}) with suitable
boundary conditions. We assume that the tangential component of $A$ vanishes on $\partial D$, i.e.
\begin{equation}\label{bcA}
n\times A=0 \quad\hbox{on }\partial D \times (0,T)\,,
\end{equation}
where $n$ denotes the outward unit normal vector to $\partial D$;
for the absolute temperature we neglect the possible radiative heat transfer between the inductor and the workpiece
assuming
\begin{equation}\label{bcteta}
\kappa \frac{\partial\teta}{\partial\nu}+\eta\teta=g\quad\hbox{on }\partial \Sigma\times (0,T)\,,
\end{equation}
where $\nu$ denotes the outward unit normal vector to $\partial \Sigma$, $\eta$ stands for an heat transfer coefficient and $g$ is a  given
boundary source.

In the following sections we will study the well-posedness  for the Cauchy system  associated to
a suitable variational formulation of (\ref{eqz}), \eqref{eqA}, \eqref{defu}, (\ref{fourier}), (\ref{def:f}--\ref{bcteta}).

%%%%%%%%%%%%%%%%%%%%%%%%%%%%%%%%%%%%%%%%%%%%%%%%%%%%%%%%%%%%%%%%%%%
%%%%%%%%%%%%%%%%%%%%%%%%% mainres %%%%%%%%%%%%%%%%%%%%%%%%%%%%%%

\parskip0mm
\section{Well-posedness results}
\label{wellposed}

In this section we introduce the functional framework, the main notation, and the assumptions on the data in order to deal with
system \eqref{eqz},  \eqref{eqA}, \eqref{defu}, \eqref{fourier}, (\ref{def:f}--\ref{bcteta}). Moreover,
we state here our main results concerning well-posedness and regularity for
a suitable weak formulation of the corresponding Cauchy problem. The proofs of these results are given in Section~\ref{proof}.
For analytical reasons our analysis is restricted to the case $\sigma=\sigma(x, z)$, $\mu=\mu(x, z)$ in \eqref{def:sigma} and \eqref{def:mu}, that is we drop the explicit temperature dependence. Since the phase fraction grows with growing temperature we still maintain the effect of temperature change on the parameter values. Moreover, we consider the standard Fourier law with $\kappa$ positive constant in \eqref{fourier}. 
%%%%%%%%%%%%%%%%%%%%%%%%%%%%%%%%%%%%%%%%%%%%%%%%%%%%%%%%%%%%%%%%%

\subsection{Notation and preliminaries} We first recall the definitions of generalized $\curl$ and $\dive$ operators as well as the related
embedding results we need in the sequel.
Referring to the Figure~\ref{domain}, let $D\subseteq \RR^3$ be a bounded 
domain. Assume
that $\overline{\Omega}\subset D$, $\overline{\Sigma}\subset D$, $\overline{\Omega}\cap\overline{\Sigma}=\emptyset$
and $\partial\Omega$, $\partial\Sigma$, $\partial D$ are of class $C^{1,1}$. Define $G=\Omega\cup\Sigma$ and
the space-time domain  $Q=\Sigma\times (0,T)$.
Let us denote with the symbol $\WW$ the vector-valued counterpart of the Banach space $W$. Let us use the symbol $\LL^q(U\times(0,t))$ for
the space $L^q(0,t; \LL^q(U))$, for every $q\in[1, +\infty)$, $t\in (0,T]$, $U$ bounded and $C^{1,1}$ domain in $\RR^3$.
Let $p\in (1,+\infty)$ and $\psi\in \LL^p(D)$,  then, we write $\curl \psi\in \LL^p(D)$ if there exists
$\xi\in \LL^p(D)$ such that
\[
\int_D\psi\curl\phi\dd x=\int_D\xi\phi\dd x\,,
\]
for all $\phi\in \CC^\infty_c(D)$. We define $\curl \psi:=\xi$, as the uniquely determined vector $\xi$.
Analogously, we write $\dive\psi\in L^p(D)$ if there exists $\zeta\in L^p(D)$ such that
\[
\int_D\psi\cdot\nabla\eta\dd x=-\int_D\zeta\eta\dd x\,,
\]
for all $\eta\in C^\infty_c(D)$. We define $\dive \psi:=\zeta$, as the uniquely determined vector $\zeta$.
Then, we introduce the Banach spaces (with the graph norms) $\LL^p_{curl}(D)=\{\psi\in \LL^p(D)\,:\,\curl\psi\in \LL^p(D)\}$ and
$\LL^p_{div}(D)=\{\psi\in \LL^p(D)\,:\,\dive\psi\in L^p(D)\}$. For $\psi\in L^2_{curl}(D)$ we define the the linear bounded
trace operator $\gamma_{\tau}(\psi)$ using the well-known Green's formula (cf., e.g. \cite{DuvautLions76})
\begin{equation}\label{green1}
\int_D\psi\curl\phi\dd x-\int_D\curl\psi\cdot\phi\dd x=-\int_{\partial D}( n\times \psi )\cdot \phi\dd s=:\duav{\gamma_{\tau}(\psi), \phi}\,,
\end{equation}
for all $\phi \in \HH^1(D)$, where (with an abuse of notation) the integral over $\partial D$ has to be understood as the duality
between $\HH^{-1/2}(\partial D)$ and $\HH^{1/2}(\partial D)$ and $n$ denotes the outward unit normal vector to $\partial D$.
Similarly, for $\psi\in L^2_{div}(D)$ we introduce the linear bounded trace operator $\gamma_n(\psi)$ by the Green's formula
\begin{equation}\label{green2}
\int_D\psi\cdot\nabla\phi\dd x+\int_D\dive\psi\cdot\phi\dd x=\int_{\partial D}( n\cdot \psi ) \phi\dd s=:\duav{\gamma_{n}(\psi), \phi}\,,
\end{equation}
for all $\phi \in H^1(D)$, where (with an abuse of notation) the integral over $\partial D$ has to be understood as the duality
between $H^{-1/2}(\partial D)$ and $H^{1/2}(\partial D)$.
Finally, we introduce the Hilbert space
\[
\XX=\{v\in \LL^2_{curl}(D)\,:\, \dive v=0\quad\hbox{and } \gamma_\tau(v)=0\}\,.
\]
Notice that, since $\partial D \in C^{1,1}$, then the space $\XX$, equipped with the norm
\[
\|v\|_{\XX}=\|\curl v\|_{\LL^2(D)}\,,
\]
is a closed subspace of $\HH^1(D)$. Moreover,
let us notice that, from the Green's formula
it follows that: if $A\in \XX$, then $\gamma_n(A)=0$. Indeed, if we denote by $w=\curl A\in \LL^2(D)$, then $\dive w=0$
and $\gamma_n(w)\in H^{-1/2}(\partial D)$. Moreover the following Green's formulas hold true for
every test function $\phi \in H^1(D)$:
\begin{align}\label{green3}
\duav{\gamma_n(w),\phi}=&\int_D\left(\dive w\phi+w\cdot\nabla\phi\right)\dd x
=\int_D\curl A\cdot\nabla\phi\dd x\\
\no
=&\int_D\curl A\cdot\nabla\phi-\int_D A\cdot\curl\nabla\phi\dd x\\
\no
=&\duav{\gamma_\tau(A), \nabla\phi}\,.
\end{align}
Hence, if $A\in \XX$, then $\gamma_n(\curl A)=0$ (cf., e.g., \cite[Chapter 2]{GiraultRaviart86}).

We recall here some results which will be useful in the sequel of the paper. The first one is an embedding result
which is a consequence of \cite[Thm. 3.3]{gries90}. The complete proof of a more general result can be found in
\cite[Prop.~2.2, p.~7]{druet07}.
\bele\label{embedd}
Let $U\subset\RR^3$ be a bounded  $C^1$ domain, then the space
\[
\mathcal{W}^{p,\alpha}(U)=\{u\in \LL^p_{curl}(U)\cap \LL^p_{div}(U)\,:\, \gamma_n(u)\in L^\alpha(\partial U)\}
\]
continuously embeds in the space $\LL^\xi(U)$ for $\xi:=\min\{\frac{3\alpha}{2}, p^*\}$, where $p^*$ is the Sobolev
embedding exponent
\[
p^*:=\begin{cases}
\frac{3p}{3-p}\quad &\hbox{if } 1\leq p<3\,,\\
s\in [1, +\infty)\quad &\hbox{if } p=3\,,\\
+\infty\quad &\hbox{if } p>3\,.
\end{cases}
\]
\enle

Moreover, we recall the following interpolation inequality, holding true for $p,q,r\in [1,+\infty]$, with $p<r<q$,
$1/r=\alpha/p+(1-\alpha)/q$, $\alpha\in (0,1)$, and $v\in \LL^q(D)$:
\begin{equation}\label{interp}
\|v\|_{\LL^r(D)}\leq \|v\|_{\LL^p(D)}^\alpha\|v\|_{\LL^q(D)}^{1-\alpha}\,.
\end{equation}

\subsection{Hypotheses} We list here our basic assumptions on the functions $\sigma$, $\mu$, $J_0$, $u$, $\tau$, and $z_{eq}$ in
\eqref{eqA}, \eqref{eqz}, and \eqref{eqteta}, where we take the constants $c_v$, $\kappa$, $\eta$, and $L$ equal to 1, for simplicity.
\behy\label{hyp1} Assume that
\begin{itemize}
\item[(i)] $\sigma(x, z)\,:\,\overline{D}\times [0,1]\to \RR$ is a continuous and Lipschitz
continuous (w.r.t. $z$ for almost all $x\in D$) function  defined in \eqref{def:sigma} and \eqref{def:sigmaz},
where $\sigma_w\in C^{1,1}([0,1])$, $\sigma_i>0$, constant.
Assume moreover that there exists two constants $0<\underline{\sigma}\leq \overline{\sigma}$ such that
\[
\underline{\sigma}\leq\sigma(x,z)\leq\overline{\sigma}\quad\hbox{in } \overline{D}\times [0,1]\,;
\]
\item[(ii)] $\mu(x, z)\,:\,\overline{D}\times [0,1]\to \RR$ is a continuous and Lipschitz continuous (w.r.t. $z$ for almost all $x\in D$)
function defined in \eqref{def:mu}, \eqref{def:muz},
where $\mu_w\in C^{1,1}([0,1])$, $\mu_i,\,\mu_0\in \RR^+$.
Assume moreover that there exists two constants $0<\underline{\mu}\leq \overline{\mu}$ such that
\[
\underline{\mu}\leq\mu(x,z)\leq\overline{\mu}\quad\hbox{in } \overline{D}\times [0,1]\,;
\]
\item[(iii)] $u\in H^1(0,T)$, $J_0\,:\, D\to \RR^3$ is an $\LL^2_{curl}(D)$-function defined in \eqref{def:J0};
\item[(iv)] $\tau, \, z_{eq}\in C^2(\RR)$ and there exists positive constants $0<\tau_*\leq \tau^*$ and $M>0$ such that
\begin{align}\no
&\tau_*\leq \tau(\teta)\leq \tau^* , \, 0\leq z_{eq}(\teta)\leq 1
\,\hbox{ for all }\teta\in \RR, \, \|\tau\|_{C^2(\RR)}\leq M,\, \|z_{eq}\|_{C^2(\RR)}\leq M\,,\\
\no
&|-z_{eq}(\teta)+\teta z_{eq}'(\teta)|\leq M \quad\hbox{for all }\teta\in \RR, \quad |\teta z_{eq}''(\teta)|\leq M\,;
\end{align}
\item[(v)] $g\in L^\infty(0,T; L^\infty(\partial\Sigma))$;
\item[(vi)] $A_0\in \XX\cap \HH^3(D)$, $\teta_0\in W^{2,5/3}(\Sigma)$.
\end{itemize}
\enhy

\beos\label{f}
We will use in the following weak formulation \eqref{eqtetaw} of the internal energy balance \eqref{eqteta}  the notation $f(\teta, z)=-L(z_{eq}(\teta)-z-\teta z'_{eq}(\teta))$. Indeed, using \eqref{eqz}  we can rewrite $-\mathcal{F}(\teta,z)z_t= -f(\teta,z)z_t$ for $z\geq z_{eq}(\teta)$ and $-\mathcal{F}(\teta,z)z_t= 0=-f(\teta,z)z_t$ for $z< z_{eq}(\teta)$. Notice moreover that   Hyp.~\ref{hyp1} (iv) implies that
\begin{equation}\label{propf}
f\in C^{0,1}(\RR\times[0,1]), \quad |f(\teta, z)|\leq C_f\quad \forall (\teta, z)\in \RR\times [0,1]\,,
\end{equation}
for some positive constant $C_f$ depending on $M$. Finally, let us observe that the assumption \eqref{propf}
are the only ones we need in order to get our next results, i.e., we do not need the explicit form \eqref{def:f}
of $f$ which was assumed in the previous section in order to comply with Thermodynamics. However,
other thermodynamically-compatible  choices are possible.
\eddos
%Check the regularity  needed on $g4 and $\teta_0$ in order to aplly maximal regularity results... Pruss e vedi anche come
%scrivere la regolarita' di teta e z

\subsection{Weak formulation and main theorem.} We are ready now to state the weak formulation of the Cauchy problem
\eqref{eqz},  \eqref{eqA}, \eqref{defu}, \eqref{fourier}, (\ref{def:f}--\ref{bcteta}).

\bepro~\label{problem} Find a triple $(A, \teta, z)$ with the regularity properties
\begin{align}
\label{regA}
& A\in H^2(0,T; \LL^2(D))\cap W^{1,\infty}(0,T;\XX), \quad \curl A\in L^\infty(0,T;\LL^6(D))\,,\\
\label{regteta}
&\teta\in W^{1,5/3}(0,T; L^{5/3}(\Sigma))\cap L^{5/3}(0,T; W^{2,5/3}(\Sigma))\cap L^2(0,T; H^1(\Sigma))\cap L^\infty(Q)\,,\\
\label{regz}
&z\in W^{1,\infty}(0,T;W^{1,^\infty}(\Sigma)),\quad 0\leq z<1\hbox{ a.e. in }Q \,,
\end{align}
solving the following system
\begin{align}
\label{eqAw}
&\int_G \sigma (x, z) A_t\cdot v\dd x+\int_D\frac{1}{\mu(x,z)}\curl A\cdot \curl v\dd x
=\int_\Omega J_0(x) u(t)\cdot v\dd x\\
\no
&\qquad\qquad\qquad\qquad\hbox{ for all }v\in \XX, \hbox{ a.e. in }(0,T)\,,\\
\label{eqtetaw}
&\teta_t-\Delta\teta=-f(\teta, z)z_t+\sigma(x, z)|A_t|^2\quad\hbox{a.e. in }Q\,,\\
\label{eqzw}
&z_t=\frac{1}{\tau(\teta)}\left(z_{eq}(\teta)-z\right)^+\quad\hbox{a.e. in }Q\,,\\
\label{bctetaw}
&\frac{\partial\teta}{\partial\nu}+\teta=g\quad \hbox{a.e. on }\partial\Sigma\times (0,T)\,,\\
\label{iniw}
&A(0)=A_0,\quad\hbox{a.e. in }D,\quad \teta(0)=\teta_0, \quad z(0)=0\quad\hbox{a.e. in }\Sigma\,.
\end{align}
and satisfying the following estimate
\begin{align}\label{estiS}
&\|A\|_{H^2(0,T; \LL^2(D))\cap W^{1,\infty}(0,T;\XX)}+\|\curl A\|_{L^\infty(0,T;\LL^6(D))}\\
\no
&+\|\teta\|_{W^{1,5/3}(0,T; L^{5/3}(\Sigma))\cap L^{5/3}(0,T; W^{2,5/3}(\Sigma))\cap L^2(0,T; H^1(\Sigma))\cap L^\infty(Q)}
\\
\no
&+\|z\|_{W^{1,\infty}(0,T;W^{1,^\infty}(\Sigma))}\leq S
\end{align}
where the constant $S$ depends on the data of the problem.
\empro
We are now ready to state our main result.
\bete\label{th1}
Let Hypothesis~\ref{hyp1} hold true. Then, there exists a unique solution to Problem~\ref{problem}. Moreover,
if we denote by $(A_i, \teta_i, z_i)$ ($i=1,2$) two triples of solutions corresponding to data $(A_{0,i}, \teta_{0,i}, u_i)$, then,
there exists a positive constant $C=C(S)$ (cf. estimate \eqref{estiS}) such that the following stability estimate holds true
\begin{align}
\label{stability}
&\|(A_1-A_2)(t)\|_{\LL^2(D)}^2+\|\curl(A_1-A_2)\|_{\LL^2(D\times (0,T))}^2\\
\no
&
\quad+\|\partial_t(A_1-A_2)(t)\|_{\LL^2(D)}^2+\|\curl(\partial_t(A_1-A_2))\|_{\LL^2(D\times (0,T))}^2\\
\no
&\quad+\|(\teta_1-\teta_2)(t)\|_{L^2(\Sigma)}^2
+\|\teta_1-\teta_2\|_{L^2(0,T; H^1(\Sigma))}^2\\
\no
&\quad+\|(z_1-z_2)(t)\|_{H^1(\Sigma)}^2+\|\partial_t(z_1-z_2)\|_{L^2(0,T; H^1(\Sigma))}^2\\
\no
&
\leq C\left( \|A_{0,1}-A_{0,2}\|_{\XX}^2+\|(\partial_t(A_1-A_2))(0)\|_{\LL^2(D)}^2+\|\teta_{0,1}-\teta_{0,2}\|_{L^2(\Sigma)}^2\right.\\
\no
&\quad
\left.+\|u_1J_0-u_2J_0\|_{L^2(0,T)}^2+\|u_1'J_0-u_2'J_0\|_{L^2(0,T)}^2  \right)\quad \hbox{for all } t\in [0,T]\,.
\end{align}
\ente

%\clearpage
%%%%%%%%%%%%%%%%%%%%%%%%%%%%%%%%%%%%%%%%%%%%%%%%%%%%%%%%%%%%%%%%%%%%%%%%%%%%%%%%%%%%%%%%%%%
%%%%%%%%%%%%%%%%%%%%%%%%%%%%%%%%%%%% proofs %%%%%%%%%%%%%%%%%%%%%%%%%%%%%%%%%%%%%%%%%%%%%
\section{Proof of Theorem~\ref{th1}}
\label{proof}

In this Section we prove Theorem~\ref{th1} in three steps: first we prove
local existence of solutions by means of Schauder fixed point argument. Secondly we
prove the global a-priori estimates necessary in order to extend the solution to
the whole time interval $[0,T]$. Finally, we prove the stability estimate \eqref{stability} entailing,
in particular, uniqueness of solution to Problem~\ref{problem}.

\subsection{Existence of a (local in time) solution}
We start the proof by solving (locally in time) Problem~\ref{problem}
by means of a standard fixed point argument of Schauder type.

For a fixed $T_0\in(0,T]$ (which shall be specified later on)
and a fixed constant $R>0$,
 let us introduce the space
\[
\OO=\{\teta\in L^2(0,T_0;L^2(\Sigma))\,:\,\|\teta\|_{L^2(0,T_0;L^2(\Sigma))}\leq R\}\,.
\]

In the following,  we shall construct an operator $\mathcal{T}$,
which
 maps $\OO$ into itself
 for a suitable  time  $0 < T_0 \leq T$,
  in such a way that any
  fixed point of $\mathcal{T}$ yields a
 solution to Problem~\ref{problem}.
  We shall prove that $\mathcal{T}$ is compact and
 continuous w.r.t. the topology of $L^2(0,T_0;L^2(\Sigma))$.
  Hence, by the Schauder theorem $\mathcal{T}$ admits (at
least) a fixed point $\teta$ in
$\OO$, whence
the existence  of a solution $(A,\teta, z)$
  to the Cauchy Problem~\ref{problem}
 on the interval~$[0,T_0]$. Finally, the (local) uniqueness result will be a consequence of the stability estimate \eqref{e15} below.

\smallskip
\noindent
{\underline{\it Definition of the fixed point map $\mathcal{T}$.}} We construct the operator $\mathcal{T}$ in this way: given
$\bar\teta\in \OO$, the operator
\[\mathcal{T}\,:\, \OO \to L^{5/3}(0,T; W^{2,5/3}(\Sigma))\cap W^{1,5/3}(0,T; L^{5/3}(\Sigma))\]
\[\mathcal{T}\,:\,  \bar \teta\mapsto\teta\quad\hbox{ where $\teta$ solves}\]
\begin{align}
\label{eqtetas}
&\teta_t-\Delta\teta=-f(\bar\teta, z)z_t+\sigma(x, z)|A_t|^2\quad\hbox{a.e. in }Q\,,\\
\label{bctetas}
&\frac{\partial\teta}{\partial\nu}+\teta=g\quad \hbox{a.e. on }\partial\Sigma\times (0,T)\,,\\
\label{initetas}
& \teta(0)=\teta_0, \quad\hbox{a.e. in }\Sigma\,,
\end{align}
and $z\in W^{1,\infty}(0,T;L^\infty(\Sigma))$, $A\in H^1(0,T;\XX)\cap W^{1,\infty}(0,T;\LL^2(D))$ solve
\begin{align}
\label{eqzs}
&z_t=\frac{1}{\tau(\bar\teta)}\left(z_{eq}(\bar\teta)-z\right)^+\quad\hbox{a.e. in }Q,\quad z(0)=0\quad\hbox{a.e. in }\Sigma\,,\\
\label{eqAs}
&\int_G \sigma (x, z) A_t\cdot v\dd x+\int_D\frac{1}{\mu(x,z)}\curl A\cdot \curl v\dd x
=\int_\Omega J_0(x) u(t)\cdot v\dd x\\
\no
&\qquad\qquad\qquad\qquad\hbox{ for all }v\in \XX, \hbox{ a.e. in }(0,T),\quad A(0)=A_0,\quad\hbox{a.e. in }D\,.
\end{align}

\smallskip
\noindent{\underline{\it $\mathcal{T}$ maps $\OO$ in itself.}} Notice that, given $\bar\teta\in \OO$,
\cite[Lemma~2.5, p.~1093]{HoembergSokolowski03} ensures that there exists
a unique solution $z\in W^{1,\infty}(0,T_0;L^\infty(\Sigma))$ to \eqref{eqzs} such that
\[
\|z\|_{W^{1,\infty}(0,T_0;L^\infty(\Sigma))}\leq C_1\,,
\]
with $C_1$ independent of $\bar\teta$, and  $0\leq z<1$ a.e. in $\Sigma$.

Then, given $z\in W^{1,\infty}(0,T_0;L^\infty(\Sigma))$, it is
possible to find a unique solution $A\in H^1(0,T_0;$ $\XX) \cap W^{1,\infty}(0,T_0;\LL^2(D))$
to \eqref{eqAs} by means of a standard implicit time-discretization scheme (cf., e.g., the monograph \cite{kacur85}). The basic estimates are
the discrete versions  of the ones which follows from taking $v=A_t$ in \eqref{eqAs} (cf. also the following {\it Second estimate}) and
taking the time derivative of \eqref{eqAs} with $v=A_t$ (cf. also the following
and {\it Third estimate}). These give the bound
\[
\|A\|_{H^1(0,T_0;\XX)\cap W^{1,\infty}(0,T_0;\LL^2(D))}\leq C_2\,,
\]
with $C_2$ independent of $\bar\teta$. Then, using the fact that $A_t\in L^\infty(0,T_0;\LL^2(D))$,
$A_t\in L^2(0,T_0;\XX)\hookrightarrow L^2(0,T_0;\LL^6(D))$ (due to
the continuous embedding of $\HH^1(D)$ into $\LL^6(D)$) and employing
the interpolation inequality \eqref{interp} with $q=6$, $p=2$ and $r=10/3$, we get $A_t\in \LL^{10/3}(D\times (0,T_0))$.
This implies that $|A_t|^2\in L^{5/3}(D\times(0,T_0))$ on the right hand side
of \eqref{eqtetas} and
\[
\|A_t\|_{\LL^{10/3}(D\times (0,T_0))}\leq C_3\,.
\]
Applying now the standard  maximal regularity results in $L^p$-spaces (cf., e.g., \cite[Thm.~3.1, Prop.~3.3]{pruss})
to (\ref{eqtetas}--\ref{initetas}), and Hyp.~\ref{hyp1} (iv) (cf. also Remark~\ref{f}), we can conclude
that there exists a unique $\teta\in L^{5/3}(0,T_0; W^{2,5/3}(\Sigma))\cap W^{1,5/3}(0,T_0; L^{5/3}(\Sigma))$ solving
(\ref{eqtetas}--\ref{initetas}) and such that
\[
\|\teta\|_{L^{5/3}(0,T_0; W^{2,5/3}(\Sigma))\cap W^{1,5/3}(0,T_0; L^{5/3}(\Sigma))}\leq C_4\,,
\]
with $C_4$ independent of $\bar\teta$.  Moreover, testing equation \eqref{eqtetas} by $\teta$, we get
\[
\|\teta\|_{L^2(0,T_0; H^1(\Sigma))\cap L^\infty(0,T_0; L^2(\Sigma))}\leq C_5\,,
\]
with $C_5$ independent of $\bar\teta$.
This implies in particular
\[
\|\teta\|_{L^2(0,t; L^2(\Sigma))}\leq t^{1/2}\|\teta\|_{L^\infty(0,T_0;L^2(\Sigma))}\leq t^{1/2} C_5\quad \hbox{for every }t\in [0,T_0)\,.
\]
Then, if we choose $T_0$ such that $T_0^{1/2}C_5\leq R$ the map {\sl $\mathcal{T}$ maps $\OO$ into itself}.

\smallskip
\noindent
{\underline{\it $\mathcal{T}$ is compact.}} The map {\sl $\mathcal{T}$ is compact} due to the compact embedding
\[
L^{5/3}(0,T_0; W^{2,5/3}(\Sigma)\cap L^2(0,T_0; H^1(\Sigma))
\cap W^{1,5/3}(0,T_0; L^{5/3}(\Sigma)) \hookrightarrow \hookrightarrow L^2(0,T_0; L^2(\Sigma))
\]
(cf.~\cite[Thm.~4,~Cor.~5]{simon}).

\smallskip
\noindent
{\underline{\it $\mathcal{T}$ is continuous.}} {\sl The continuity of $\mathcal{T}$} follows from these facts:
\begin{itemize}
\item[1.] From \cite[Lemma~2.5]{HoembergSokolowski03}
we have that,
if $\bar\teta_k\to \teta$ strongly in $L^2(0,T_0; L^2(\Sigma))$, for $k\to +\infty$,
then $z_k\to z$ strongly in $W^{1,p}(0,T_0; L^p(\Sigma))$ for every $p\in [1,+\infty)$, where $z_k$ (and $z$) are the solutions to \eqref{eqzs}
corresponding to $\bar\teta_k$ (and $\bar\teta$), respectively.
\item[2.] From the compact embedding of
$H^1(0,T_0;\XX)\cap W^{1,\infty}(0,T_0;\LL^2(D))$
into $\LL^q(D\times (0, T_0))$ (for $q\in [1, 10/3)$), we get that $\partial_t A_k \to \partial_t A$ strongly in
$\LL^q(D\times (0, T_0))$ for every $q\in [1, 10/3)$,
and $\curl A_k\to \curl A$ weakly in $H^1(0, T_0; \LL^2(D))$ (at least for a subsequence of $k$),
where $A_k$  are the solutions to \eqref{eqAs} corresponding to $z_k$. We need now to show that the limit  $A$ is the solution of \eqref{eqAs}
corresponding to $z$. Then, the convergence will hold true for the whole sequence $A_k$.
 Indeed, from the strong convergence of $z_k$ and using Hyp.~\eqref{hyp1} (i), (ii), we get that
$\sigma(x, z_k)\to \sigma (x, z)$  and $\mu(x, z_k)\to \mu(x, z)$ strongly in $L^p(0,T_0; L^p(\Sigma))$ for every $p\in [1,+\infty)$.
This implies that we can pass to the limit for $k\to+\infty$ and the limit $A$ solves \eqref{eqAs}.
\item[3.] Finally, we can come back to equation
\eqref{eqtetas} and we have that $|\partial_tA_k|^2 \to |\partial_t A|^2$ in $L^s(\Sigma\times (0,T_0))$ for some $s\in [1, 5/4)$ because
$\partial_t A_k\to \partial_t A$ strongly in $\LL^q(D\times (0, T_0))$ for every $q\in [1, 10/3)$
and $\partial_t A_k$, $\partial_t A$ are bounded (uniformly in $k$) in $\LL^2(D\times (0, T_0))$. This fact, together with the strong convergence
$\partial_t z_k\to \partial_t z$ strongly in $L^p(0,T_0; L^p(\Sigma))$ for every $p\in [1,+\infty)$, $f(\bar\teta_k, z_k)\to f(\bar\teta, z)$
strongly in $L^2(0,T_0; L^2(\Sigma))$ (cf. Hyp.~\ref{hyp1} (iv) and Remark~\ref{f}), imply that also
 $\teta_k\to \teta$ strongly in $L^2(\Sigma\times (0,T_0))$, for $k\to +\infty$, where $\teta_k$ (and $\teta$) are
are the solutions to \eqref{eqtetas}
corresponding to $z_k$, $A_k$ (and $z$, $A$), respectively.
\end{itemize}
Hence, by the
Schauder theorem {\sl $\mathcal{T}$ admits (at least) a fixed point} $\teta$ in $\OO$, whence the existence of
a solution $(A, \teta, z)$ to the Cauchy Problem~\ref{problem} on the interval $[0,T_0]$.

\subsection{Global a-priori estimates} In order to extend the solution to the whole time interval $[0,T]$ we need to
prove suitable global (independent of $T_0$) a-priori estimates. In what follows the positive constants are denoted by the same symbol $C$
even if they are different from line to line. They may depend on $\sigma$, $\mu$, $\tau$, $z_{eq}$, $A_0$, $\teta_0$, $g$, $u$, $J_0$,
and $T$, but not on $T_0$.
\smallskip

\noindent
\underline{{\it First estimate.}} First of all, we can apply \cite[Lemma~2.5, p. 1093]{HoembergSokolowski03} to the equation \eqref{eqzw} entailing
the following estimate:
\begin{equation}\label{e1}
\|z\|_{W^{1,\infty}(0,T_0;L^\infty(\Sigma))}\leq C\,,\quad 0\leq z(x,t)< 1 \quad\hbox{a.e. in }Q\,.
\end{equation}

\smallskip
\noindent
\underline{{\it Second estimate.}} Take $t\in (0,T_0)$ and $v=A_t$ in \eqref{eqAw}. Using Hyp~\ref{hyp1} (i), (ii)
(cf.~also \eqref{def:muz}), we get (using the symbol $s$ for the time variable and $x$
for the space variable inside the integrals, but making the dependence on $x$ and $s$ explicit only when it is necessary):
\begin{align}\label{e2}
\underline{\sigma}\int_0^t\int_G |A_s|^2\dd x \dd s+&\int_0^t\int_D\frac{1}{\mu(x,z)}\curl A\cdot \curl A_s\dd x \dd s\\
\no
=&\int_0^t\int_\Omega  u(s) J_0(x)\cdot  A_s\dd x\dd s\,.
\end{align}
Now, integrating by parts in time the second term in \eqref{e2}, we get
\begin{align}\no
\int_0^t\int_D\frac{1}{\mu(x,z)}\curl A\cdot \curl A_s\dd x \dd s=&\frac12\int_0^t\int_\Sigma\frac{\mu'_{w}(z)z_s}{\mu^2(x,z)}|\curl A|^2\dd x\dd s\\
\no
&+\int_D\frac{1}{2\mu(x,z(t))}|\curl A|^2(t)\dd x\\
\no
&-\frac12\int_D\frac{1}{\mu(x,0)}|\curl A_0|^2\dd x\,.
\end{align}
Then, using estimate \eqref{e1} and Hyp.~\ref{hyp1} (ii), we get
\begin{align}\no
\frac{\underline{\sigma}}{2}\int_0^t\int_G |A_s|^2\dd x \dd s&+
\frac{1}{2\overline{\mu}}\int_D|\curl A|^2(t)\dd x\leq\frac{1}{2\underline\mu}\int_D|\curl A_0|^2\dd x\\
\no
&+C\left(\int_0^t\int_D|\curl A|^2\dd x\dd s+\itt\int_\Omega|u(s)|^2|J_0|^2\dd x\dd s\right)\,.
\end{align}
Using now a standard Gronwall lemma together with Hyp~\ref{hyp1} (vi), we obtain
\begin{equation}\label{e3}
\int_0^t\int_G |A_s|^2\dd x \dd s+\int_D|\curl A|^2(t)\dd x\leq C\left(1+\itt\int_\Omega|u(s)|^2|J_0(x)|^2\dd x\dd s\right)\,.
\end{equation}

\smallskip
\noindent
\underline{{\it Third estimate.}} We can now formally (in order to make it rigorous one should perform the estimate, e.g.,
on a time discrete scheme) differentiate \eqref{eqAw} with respect to $t$ and take $v=A_t$ as test function. Integrating over $(0,t)$, using
Hyp~\ref{hyp1} (i), (ii), as well-as the integration by parts in time in the second summand, gives (cf. also \eqref{def:sigmaz}, \eqref{def:muz}):
\begin{align}\no
\frac{\underline{\sigma}}{2}\int_G|A_t|^2(t)\dd x&+\frac{1}{2\overline{\mu}}\itt\int_D|\curl A_s|^2\dd x\dd s\leq
\overline{\sigma}\int_D|A_t(0)|^2\dd x\\
\no
&-\frac{1}{2}\itt\int_\Sigma\sigma_{w}'(z)z_s|A_s|^2\dd x \dd s+\itt\int_\Sigma\frac{\mu_{w}'(z)z_s}{\mu^2(x,z)}\curl A\cdot\curl A_s\dd x \dd s\\
\no
&+\itt\int_\Omega|u'(s)|^2|J_0(x)|^2\dd x\dd s\,.
\end{align}
Using now estimates \eqref{e1}, \eqref{e3},  and Hyp.~\ref{hyp1} (i), (ii), (vi), we obtain
\begin{align}\label{e4}
\int_G |A_t|^2(t)\dd x &+\itt\int_D|\curl A_s|^2\dd x\dd s\\
\no
&\leq C\left(1+\itt\int_\Omega\left(|u'(s)|^2+|u(s)|^2+|u(0)|^2\right)|J_0(x)|^2\dd x\dd s\right)\,,
\end{align}
where we have used the inequality $\|A_t(0)\|^2_{\LL^2(D)}\leq C(\|A_0\|_{\HH^2(D)}^2+\int_\Omega|u(0)|^2|J_0(x)|^2\dd x)$.

\smallskip
\noindent
\underline{{\it Fourth estimate.}} Collecting \eqref{e3} and \eqref{e4} and using the continuous embedding of $\HH^1(D)$ into $\LL^6(D)$) and employing
the interpolation inequality \eqref{interp} with $q=6$, $p=2$ and $r=10/3$, and Hyp.~\ref{hyp1} (iii), we get
\begin{equation}
\label{e5}
\|A\|_{H^1(0,T_0;\XX)\cap W^{1,\infty}(0,T_0;\LL^2(D))}+\|A_t\|_{\LL^{10/3}(D\times (0,T_0))}\leq C\,.
\end{equation}
This implies that $|A_t|^2\in L^{5/3}(Q)$ in the right hand side
of \eqref{eqtetaw}. Applying now the standard  maximal regularity results in $L^q$-spaces (cf., e.g., \cite[Thm.~3.1, Prop.~3.3]{pruss})
to (\ref{eqtetaw}), (\ref{bctetaw}--\ref{iniw}$_2$), together with Hyp.~\ref{hyp1} (iv), (v), (vi) we can deduce the estimate
\begin{equation}
\label{e6}
\|\teta\|_{L^{5/3}(0,T_0; W^{2,5/3}(\Sigma))\cap W^{1,5/3}(0,T_0; L^{5/3}(\Sigma))}\leq C\,.
\end{equation}
Moreover, we can test \eqref{eqtetaw} by $\teta$ obtaining, in particular, the estimate
\begin{equation}
\label{e7}
\|\teta\|_{L^{2}(0,T_0; H^1(\Sigma))\cap L^\infty(0,T_0; L^{2}(\Sigma))}\leq C\,.
\end{equation}
Moreover, from the non-negativity of the r.h.s in \eqref{eqtetaw} and of the initial and boundary conditions
we also get $\teta\geq 0$ a.e.
More in general, we also get from the standard  maximal regularity results in $L^q$-spaces (cf., e.g., \cite[Thm.~3.1, Prop.~3.3]{pruss})
\begin{equation}
\label{e8}
\|\teta\|_{L^{\frac{2p}{3(p-2)}}(0,T_0; W^{2,\frac{p}{2}}(\Sigma))\cap W^{1,\frac{2p}{3(p-2)}}(0,T_0; L^{\frac{p}{2}}(\Sigma))}\leq C\,,
\end{equation}
for all $p\in (2, 6]$.

\smallskip
\noindent
\underline{{\it Sixth estimate.}} By comparison in the ODE \eqref{eqzw}, using \eqref{e8} and Hyp.~\ref{hyp1} (iv), we get, for all
$p\in (2,6]$
\begin{equation}
\label{e8ter}
\|z_t\|_{L^{2}(0,T_0; W^{2,p/2}(\Sigma))}+\|z_{tt}\|_{L^{\frac{2p}{3(p-2)}}(0,T_0;L^{p/2}(\Sigma))}\leq C\,,
\end{equation}
implying
\begin{equation}
\no
\|\nabla z\|_{L^{\infty}(0,T_0; W^{1,p/2}(\Sigma))}\leq C\,,
\end{equation}
and so
\begin{equation}
\label{e9}
\|\nabla z\|_{L^{\infty}(\Sigma\times (0,T_0))}\leq C\,.
\end{equation}

\smallskip
\noindent
\underline{{\it Sixth estimate.}} Using the $L^\infty(0,T_0; \LL^2(D))$-regularity of $A_t$ (cf.~\eqref{e5}) and Hyp.~\ref{hyp1} (iii)
on $u$ and $J_0$, by comparison in \eqref{eqAw}, we get
\[
\left\|\curl\left(\frac{1}{\mu(x,z)}\curl A\right)\right\|_{L^\infty (0,T_0; \LL^2(D))}\leq C\,.
\]
Moreover, we can formally compute the divergence operator of $\frac{1}{\mu(x,z)}\curl A$ and we get (applying Hyp.~\ref{hyp1} (ii))
\[
\dive\left(\frac{1}{\mu(x,z)}\curl A\right)=\frac{\mu_w'(z)\nabla z}{\mu^2(x,z)}\curl A\in L^\infty(0,T_0; \LL^{2}(\Sigma))\,,
\]
because $\curl A\in L^\infty(0,T_0; \LL^2(D))$ and $\nabla z\in L^\infty(0,T; L^\infty(\Sigma))$ due to \eqref{e9}.
Thus, applying Lemma~\ref{embedd} with $p=2$, $U=D$, $\xi=6$, we obtain
\begin{equation}\label{e10}
\|\curl A\|_{L^\infty(0,T_0; \LL^6(D))}\leq C\,.
\end{equation}
Notice that we can apply Lemma~\ref{embedd} with every exponent $\alpha$ (hence $\xi=p^*=6$) because, due to formula \eqref{green3} and
to the fact that $A\in L^\infty (0,T;\XX)$, we have that $\gamma_n(\curl A)=0$.

\smallskip
\noindent
\underline{{\it Seventh estimate.}} We can now formally (in order to make it rigorous one should perform the estimate, e.g.,
on a time discrete scheme) differentiate \eqref{eqAw} with respect to $t$ and take $v=A_{tt}$ as test function. Integrating over $(0,t)$
and using Hyp~\ref{hyp1} (i), (ii) (cf.~\eqref{def:sigmaz}, \eqref{def:muz})
we get
\begin{align}\no
\underline{\sigma}\itt\int_D|A_{ss}|^2\dd x\dd s+\frac{1}{2\overline{\mu}}\int_D|\curl A_t|^2(t)\dd x
=&\frac{1}{2\underline{\mu}}\int_D|\curl A_t(0)|^2\dd x\\
\no
&-\itt \int_\Sigma\frac{\mu_w'(z)}{\mu^2(x, z)} z_s|\curl A_s|^2\dd x\dd s\\
\no
&-\itt\int_\Sigma \sigma_w'(z) z_s A_s A_{ss}\dd x\dd s\\
\no
&-\itt\int_\Sigma \frac{\mu_w'(z)}{\mu^2(x, z)} z_s\curl A_{ss}\curl A\dd x\dd s\\
\no
&+\itt \int_\Omega J_0(x) u'(s)\dd x\dd s\,.
\end{align}
We start estimating the second integral in the r.h.s using Hyp~\ref{hyp1} (ii)  and estimates \eqref{e1}, \eqref{e5} as follows:
\[
-\itt \int_\Sigma\frac{\mu_w'(z)}{\mu^2(x, z)} z_s|\curl A_s|^2\dd x\dd s\leq C \|\curl A_t\|_{L^2(0,T_0; \LL^2(D))}\leq C\,.
\]
In order to bound the third integral we use Hyp.~\ref{hyp1} (i) and again estimates \eqref{e1}, \eqref{e5} as follows:
\begin{align}\no
-\itt\int_\Sigma \sigma_w'(z) z_s A_s A_{ss}\dd x\dd s&\leq C\|A_t\|_{L^2(0,T_0; \LL^2(D))}^2
+\frac{1}{2\underline{\sigma}}\itt\int_D|A_{ss}|^2\dd x\dd s\\
\no
&\leq C+\frac{\underline{\sigma}}{2}\itt\int_D|A_{ss}|^2\dd x\dd s\,.
\end{align}
In the fourth integral we need to use the integration by parts formula, Hyp.~\ref{hyp1} (ii), and estimates \eqref{e1}, \eqref{e5},
and \eqref{e10}:
\begin{align}\no
&\itt\int_\Sigma \frac{\mu_w'(z)}{\mu^2(x, z)} z_s\curl A_{ss}\curl A\dd x\dd s\\
\no
&\leq C
\itt \|\curl A_s\|_{\LL^2(D)}\|\curl A\|_{\LL^6(D)}\|z_{ss}\|_{L^3(\Sigma)}\dd s+C\itt\|\curl A_s\|_{\LL^2(D)}^2\dd s
\\
\no
&\leq C\itt \|\curl A_s\|_{\LL^2(D)}\|z_{ss}\|_{L^3(\Sigma)}\dd s+C\itt\|\curl A_s\|_{\LL^2(D)}^2\dd s\,.
\end{align}
Now, we use Hyp.~\ref{hyp1} (iii), (vi), estimate \eqref{e8ter} with $p=3$, and apply a standard Gronwall lemma to the following
inequality
\[
\frac{\underline{\sigma}}{2}\itt\int_D|A_{ss}|^2\dd x\dd s+\frac{1}{2\overline{\mu}}\int_D|\curl A_t|^2(t)\dd x
\leq C+C\itt \|\curl A_s\|_{\LL^2(D)}\|z_{ss}\|_{L^3(\Sigma)}\dd s\,,
\]
getting the desired estimate
\begin{equation}\label{e11}
\|A_{tt}\|_{L^2(0,T_0;\LL^2(D))}+\|\curl A_t\|_{L^\infty(0,T_0;\LL^2(D))}\leq C\,.
\end{equation}
\smallskip
and so, since $|A_t|^2\in L^\infty(0,T_0; \LL^3(D))$, using \cite{kawanago},
we get
\begin{equation}
\label{e11bis}
\|\teta\|_{L^{\infty}(\Sigma\times (0,T_0))}\leq C\,.
\end{equation}
\noindent
Collecting estimates \eqref{e1}, \eqref{e5}, \eqref{e6}, \eqref{e10}, \eqref{e11}, and \eqref{e11bis},
we can now extend the solution we found on $[0,T_0]$ to
the whole time interval $[0,T]$. Finally, notice that testing \eqref{eqtetaw} by $-\teta^{-}$ ($(\cdot)^-$ denoting the negative part),
and using the definition of $f$ \eqref{def:f} together with Hyp.~\ref{hyp1} (v), (vi), we gain the non-negativity of $\teta$.
This concludes the proof of existence of solutions to Problem~\ref{problem}.

\subsection{Stability estimate} In this part we prove the stability estimate \eqref{stability}, entailing,
in particular, uniqueness of solutions to Problem~\ref{problem}.
Consider the time derivative of equation \eqref{eqAw} and rewrite it in the following form
\begin{align}\label{eqAwt}
&\int_G \left(h(x, z, z_t)A_t+\sigma(x, z) A_{tt}\right) \cdot v\dd x
\\
\no
&\qquad+\int_D\left(\frac{1}{\mu(x,z)}\curl A_t+\ell(x, z, z_t)\curl A\right)\cdot \curl v\dd x=\int_\Omega J_0(x) u'(t)\cdot v\dd x\,,
\end{align}
for all $v\in \XX$ and a.e. in $(0,T)$, where
$$
h(x, z, z_t):= \sigma'(x,z) z_t, \quad \ell(x, z, z_t):=\frac{\mu_z(x,z) z_t}{\mu^2(x, z)}\,.
$$
Let
$(A_i, \teta_i, z_i)$ ($i=1,2$) be two triples of solutions corresponding to data $(A_{0,i}, \teta_{0,i}, u_i)$.
Take the difference between \eqref{eqAwt}$_1$ and \eqref{eqAwt}$_2$ and
take the test function $v=\bar{A}_t:=(A_1-A_2)_t$. Then, denoting by $\bar{z}:=z_1-z_2$, we get, integrating over $(0,t)$
\begin{align}\label{e12}
&\frac{\underline{\sigma}}{2}\int_G|\bar{A}_t(t)|^2\, \dd x+\frac{1}{\overline{\mu}}\itt\int_D|\curl \bar{A}_s|^2\, \dd x\\
\no
&\leq \frac{\underline{\sigma}}{2}\int_G|\bar{A}_t(0)|^2\, \dd x
+ C_\sigma\itt\left(\|\bar{z}\|_{L^\infty(\Sigma)}+\|\bar{z}_s\|_{L^\infty(\Sigma)}\right)\|\bar{A}_s\|_{\LL^2(D)}^2\, \dd s\\
\no
&\qquad+C_\sigma \itt\left(\|\bar{A}_s\|_{\LL^2(D)}^2+\|\bar{z}\|_{ L^\infty(\Sigma)}\|A_{1,ss}\|_{\LL^2(D)}\|\bar{A}_s\|_{L^2(D)}\right)\, \dd s\\
\no
&\qquad+C_\mu \itt \left(\|\bar{z}\|_{L^\infty(\Sigma)}+\|\bar{z}_s\|_{L^\infty(\Sigma)}\right)\|\bar{A}_s\|_{\LL^2(D)}^2\, \dd s\\
\no
&\qquad+C_\mu\itt\left(\|\curl\bar{A}\|_{\LL^2(D)}\|\bar{A}_s\|_{\LL^2(D)}+\|\bar{z}\|_{ L^\infty(\Sigma)}\|\curl\bar{A}_s\|_{\LL^2(D)}\right)\, \dd s\\
\no
&\qquad+C\itt \|u_1'-u_2'\|\bar{A}_s\|_{\LL^2(D)}\,\dd s,
\end{align}
where we have used the following inequalities (holding true for a.e. $x$) due to Hyp.~\ref{hyp1} (i), (ii):
\begin{align}\no
&|h(x, z_1, z_{1,t})- h(x, z_2, z_{2,t})|\leq C_\sigma \left(|z_{1,t}||z_1-z_2|+|z_{1,t}-z_{2,t}|\right)\,,\\
\no
&|h(x, z_2, z_{2,t}|\leq C_\sigma |z_{2, t}|\,,\\
\no
&|\ell(x,z_1, z_{1,t})-\ell(x, z_2, z_{2,t})|\leq C_\mu \left(|z_1-z_2||z_{1,t}|+|(z_1-z_2)_t|\right)\,,\\
\no
&|\ell(x, z_2, z_{2,t}|\leq C_\mu |z_{2, t}|\,.
\end{align}
Test now \eqref{eqtetaw}$_1$-\eqref{eqtetaw}$_2$ by $\bar\teta:=\teta_1-\teta_2$, integrate over $(0,t)$,  and  use the following inequalities
\begin{align}\label{e13}
\|\bar z\|_{ L^\infty(\Sigma)}&+\itt \|\bar{z}_s\|_{L^\infty(\Sigma)}\, \dd s\leq C\itt \|\bar\teta\|_{L^\infty(\Sigma)}\, \dd s\leq C\itt \|\bar\teta\|_{W^{2,3/2}(\Sigma)}\, \dd s\\
\no
&\leq C \itt \left(\||A_{1,s}|^2-|A_{2,s}|^2\|_{\LL^{3/2}(D)}+\|\bar{z}_s\|_{L^2(\Sigma)}\right)\, \dd s\\
\no
&\leq C\itt \left(\left(\|A_{1,s}\|_{\LL^6(D)}+\|A_{2,s}\|_{\LL^6(D)}\right)\|\bar{A}_s\|_{\LL^2(D)}+\|\bar\teta\|_{L^2(\Sigma)}\right)\, \dd s\,,
\end{align}
obtaining (for all positive $\delta$)
\begin{align}\label{e14}
&\|\bar{\teta}(t)\|_{L^2(\Sigma)}^2+\itt\|\nabla\bar{\teta}\|_{L^2(\Sigma)}^2\, \dd s\leq C\|\bar{\teta}(0)\|_{L^2(\Sigma)}^2\\
\no
&\quad +\delta \itt \|\bar{\teta}\|_{L^6(\Sigma)}^2\, \dd s+C_\delta \itt \left(\|\bar{\teta}\|_{L^2(\Sigma)}^2
+\left(\|A_{1,s}\|_{\LL^6(D)}+\|A_{2,s}\|_{\LL^6(D)}\right)\|\bar{A}_s\|_{\LL^2(D)}\right)\, \dd s\,.
\end{align}
Add now the term $\itt \|\bar{\teta}\|_{L^2(\Sigma)}^2$ to both sides of \eqref{e14}, choose
$\delta$ sufficiently small and sum the result  up to \eqref{e12}, getting (for all positive constant $\eta$)
\begin{align}\label{e15}
&\|\bar{A}_t(t)\|_{\LL^2(D)}^2+\itt \|\curl\bar{A}_t\|_{\LL^2(D)}^2\, \dd s+ \|\bar\teta(t)\|_{L^2(\Sigma)}^2+\itt \|\bar\teta\|_{H^1(\Sigma)}^2\, \dd s\\
\no
&\leq C\|\bar{A}_t(0)\|_{\LL^2(D)}^2+C\|\curl \bar {A}(0)\|_{\LL^2(D)}^2+C\|\bar\teta(0)\|_{L^2(\Sigma)}^2+C \itt |u_1(s)-u_2(s)|^2\, \dd s\\
\no
&\qquad+\eta\itt \|\curl \bar{A}_s\|_{\LL^2(D)}^2\, \dd s+C_\eta\itt \|\bar\teta\|_{L^2(\Sigma)}^2\, \dd s
\\
\no
&\qquad+C_\eta\itt \left(\|A_{1,s}\|_{\LL^6(D)}^2+\|A_{2,s}\|_{\LL^6(D)}^2+\|A_{1,ss}\|_{\LL^2(D)}^2\right)\|\bar{A}_s\|_{\LL^2(D)}^2\, \dd s\,,
\end{align}
where we have used once more \eqref{e13}.
Choosing now $\eta$ sufficiently small in \eqref{e15} and applying a standard Gronwall lemma
together with the regularity properties \eqref{regA} of the solution, we obtain the desired stability estimate \eqref{stability}.
This concludes the proof of Theorem~\ref{th1}.

%%%%%%%%%%%%%%%%%%%%%%%%%%%%%%%%%%%%%%%%%%%%%%%%%%%%%%%%%%%%%%%%%%%%%%%%%%%%%%%%%%%%%%%%%%%%%%%%%%%%%%

\section{Numerical simulations}
\label{simulations}

In this section, numerical simulations for Problem \ref{problem} are carried out. We use the same setting as in Section \ref{wellposed}, all material parameters are assumed to depend on the phase fraction $z$ but not on the temperature $\teta$. The physical, temperature dependent parameters taken from literature are replaced by average values over the typical temperature range for induction hardening processes.
For further numerical simulations and an experimental verification for the complete temperature dependent induction hardening problem where also mechanical strains and stresses are considered, we refer to \cite{mefresim:2014}.

The coupled system consisting of the vector potential equation, the heat equation and the rate law for the phase fraction is solved using the finite element method. For this, we use the software {\tt pdelib} that is developed at WIAS. In order to discretize the temperature $\teta$ we use standard$P1$ elements. As  space for the vector potential $A$ we take the Hilbert space $H(\text{curl},D) = \LL^2_{curl}(D)$. In order to discretize $A$, curl-conforming finite elements of N\'ed\'elec type are used \cite{nedelec:80}. These admit only tangential continuity across interelement boundaries while the normal component, especially at material interfaces, might be discontinuous. Instead of nodal values as for classical finite elements, the degrees of freedom are the tangential components along edges of the underlying triangulation of the domain $D$.

Regarding the discretization in space, one has to account for the skin effect: alternating currents tend to distribute in a small surface layer of the workpiece. The computational grid has to resolve this small surface region. This is resolved by creating an adaptive grid, that has a high resolution in the surface area of the workpiece and is coarse in regions where a high accuracy is not necessary.

Regarding the time discretization, the system has to be solved using different time steps, since the vector potential equation and the heat equation admit different time scales.
If we consider the vector potential equation for a given temperature $\teta$, then it represents a parabolic equation that admits a time periodic solution due to the time periodic source term. We solve this equation for some periods using an order two time stepping scheme with time step $\delta t$.

\begin{figure}[!h]
\centering
\includegraphics[trim=0mm 0mm 35mm 15mm, clip, width=0.26\linewidth]{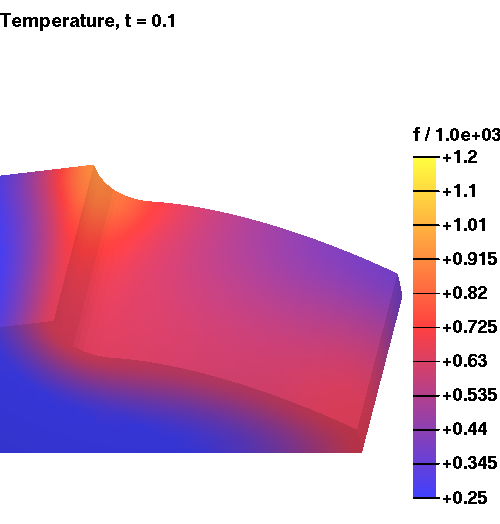}\hspace{2mm}%
\includegraphics[trim=0mm 0mm 35mm 15mm, clip, width=0.26\linewidth]{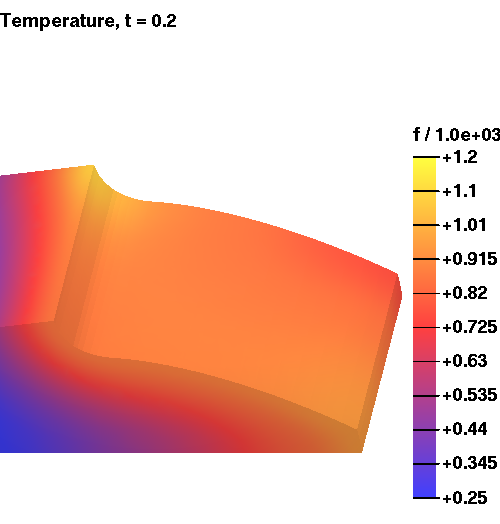}\hspace{2mm}%
\includegraphics[trim=0mm 0mm 35mm 15mm, clip, width=0.26\linewidth]{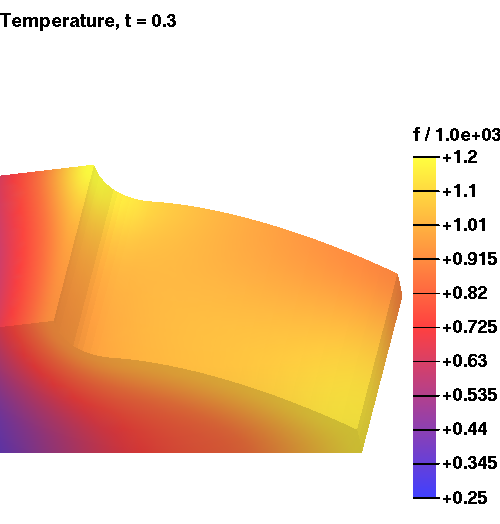}\hspace{2mm}%
\includegraphics[trim=145mm 0 0 15mm, clip, width=0.06\linewidth]{efd_z_53mf_22hf_temp-00149}\\
\subfigure[$t=0.1\,$s]{\includegraphics[trim=0mm 0mm 35mm 15mm, clip,width=0.26\linewidth]{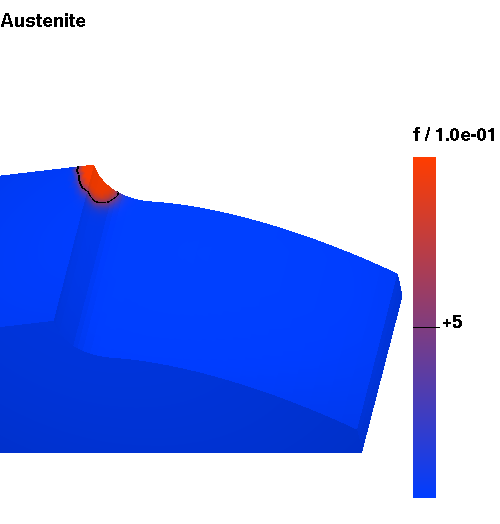}}\hspace{2mm}%
\subfigure[$t=0.2\,$s]{\includegraphics[trim=0mm 0mm 35mm 15mm, clip,width=0.26\linewidth]{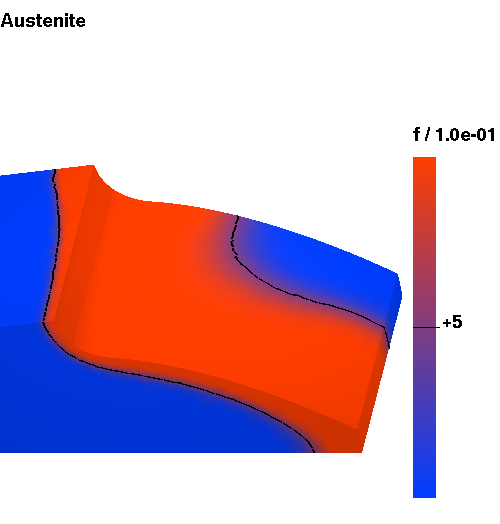}}\hspace{2mm}%
\subfigure[$t=0.3\,$s]{\includegraphics[trim=0mm 0mm 35mm 15mm, clip,width=0.26\linewidth]{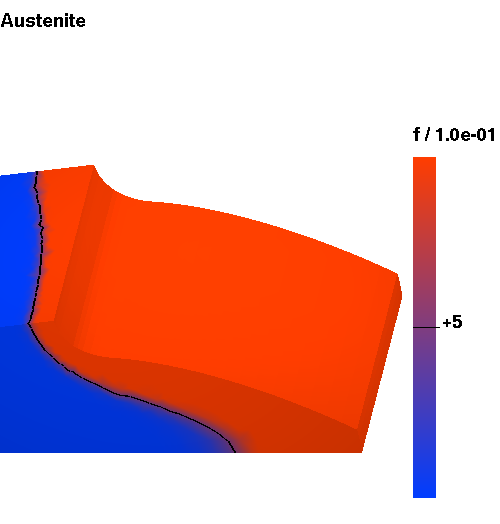}}\hspace{2mm}%
\includegraphics[trim=145mm 0 0 15mm, clip, width=0.06\linewidth]{efd_z_53mf_22hf_z-00149}\\
\caption{Numerical simulation of induction heating of a gear using the
  multifrequency approach, temperature (top row) and volume fraction
  of austenite (bottom row).}
\label{fig:sim_multi}
\end{figure}

The heat equation also represents a parabolic equation with rapidly varying right hand side $\sigma|A_t|^2$. Since heat conduction is supposed to happen on a time scale that is much slower than the oscillating current, the temperature changes at a time scale that is much larger than that of the right hand side, which is governed by the frequency of the source current. The usual approach is to approximate the Joule heat term by its average over one period, \cite{clain:1993}. Then, the heat equation together with  the rate law describing the phase transition can be solved using time steps $\Delta t \gg \delta t$, where the  rapidly varying Joule heat is replaced by an averaged Joule heat term that is obtained from the solution of the vector potential equation.

In order to reduce computational time, we make use of symmetry conditions. The cross section through the gear represents a symmetry plane, it is only necessary to consider the upper half of the geometry shown in Figure \ref{domain}. In addition, it is sufficient to perform the computations only for half of a tooth. The boundary conditions on the cutting planes are adapted accordingly.

Figure \ref{fig:sim_multi} shows the temperature and the phase
fraction of austenite for different time snapshots for a simulation
using the multifrequency approach. Starting from room temperature, the
temperature increases in the surface layer of the tooth and finally
rises above the austenitization temperature. The formation of the high
temperature phase austenite starts. After $0.3\,$s, one obtains an
austenitic profile that follows the contour of the tooth. The
corresponding temperature profile is shown in Figure
\ref{fig:sim_multi} (c), the temperature is not homogeneous across the
tooth, the highest temperatures are attained at the top edge of the
root and at the tip of the tooth. This is an important for engineers planning a heat treatment in order to prevent a melting of the gear
during the hardening process.

\begin{figure}[h]
\centering
\subfigure[MF]{\includegraphics[trim=0mm 0mm 35mm 15mm, clip,width=0.26\linewidth]{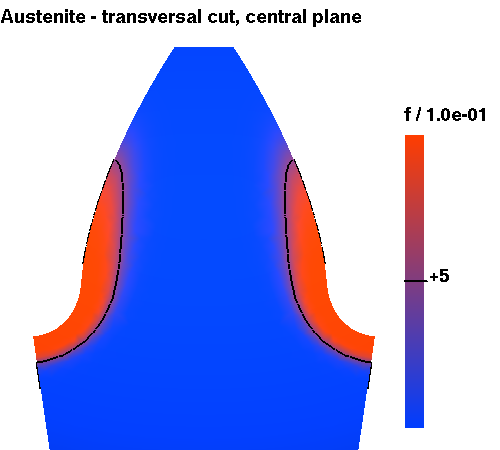}}\hspace{2mm}%
\subfigure[HF]{\includegraphics[trim=0mm 0mm 35mm 15mm, clip,width=0.26\linewidth]{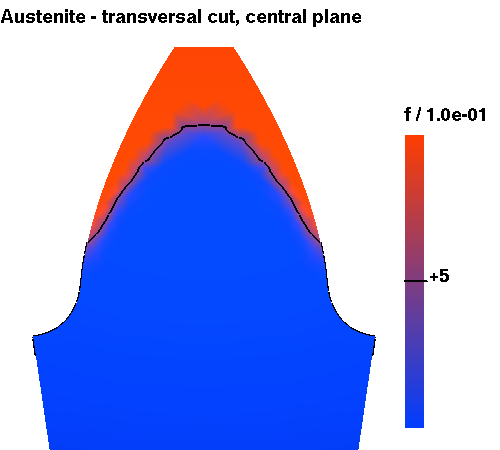}}\hspace{2mm}%
\subfigure[MF+HF]{\includegraphics[trim=0mm 0mm 35mm 15mm, clip,width=0.26\linewidth]{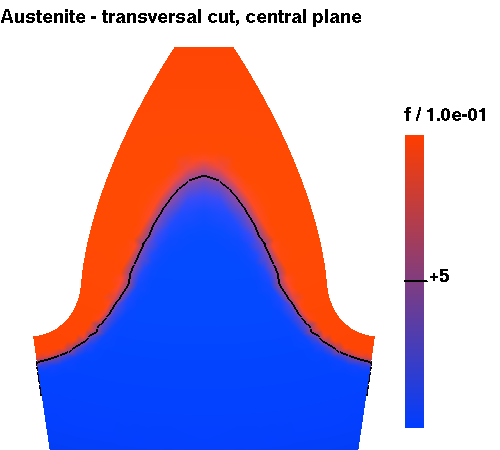}}\hspace{2mm}%
\includegraphics[trim=140mm 0 0 15mm, clip, width=0.06\linewidth]{sim_hf_cross}
\caption{Austenitic profile $z$ of a gear heated with medium-, high-
  and multifrequency power.} 
\label{fig:sim_comp}
\end{figure}

In Figure \ref{fig:sim_comp}, simulations of the austenitic fraction using only a single frequency, either MF or HF, are compared to the multifrequency case. As one can see, in the case of MF only the root of the tooth is hardened while high frequency leads to a hardening of the tip of the tooth. Only a suitable combination of MF power, HF power and heating time leads to a contour hardening of the tooth using the multifrequency approach.
The numerical computations reflect the physical aspects as explained in the introduction, see also Figure \ref{fig:multifrequency}.

%%%%%%%%%%%%%%%%%%%%%%%%%%%%%%%%%%%%%%%%%%%%% biblio %%%%%%%%%%%%%%%%%%%%%%%%%%%%%%%%%%%%%%%%%%%%%%%%%
%%%%%%%%%%%%%%%%%%%%%%%%%%%%%%%%%%%%%%%%%%%%%%%%%%%%%%%%%%%%%%%%%%%%%%%%%%%%%%%%%%%%%%%%%%%%%%%%%%%%%%

%\newpage

\end{document}